\documentclass[12pt,reqno]{amsart}

\usepackage{enumitem}   
\usepackage{graphicx}
\usepackage[colorlinks = true,
linkcolor = blue,
urlcolor  = blue,
citecolor = blue,
anchorcolor = blue]{hyperref}

\newtheorem{theorem}{Theorem}

\newtheorem{proposition}
{Proposition}
\newtheorem{corollary}
{Corollary}
\newtheorem{remark}{Remark}

\newtheorem{lemma}{Lemma}

\newfont{\bb}{msbm10 at 12pt}

\def\P{\hbox{\bb P}}

\def\SB{\mathbf{S}M}
\def\TSB{\mathcal{S}}

\def\E{\mathcal{E}}
\def\P{\mathcal{P}}

\def\<{\langle}     
\def\>{\rangle}     
\def\div{{\rm div}}

\def\Tr{{\rm tr}}

\newcommand{\bal}{\begin{align}}      \newcommand{\eal}{\end{align}}
\newcommand{\ba}{\begin{array}}      \newcommand{\ea}{\end{array}}
\newcommand{\bc}{\begin{center}}     \newcommand{\ec}{\end{center}}
\newcommand{\be}{\begin{enumerate}}  \newcommand{\ee}{\end{enumerate}}
\newcommand{\beq}{\begin{eqnarray}}  \newcommand{\eeq}{\end{eqnarray}}
\newcommand{\beQ}{\begin{eqnarray*}} \newcommand{\eeQ}{\end{eqnarray*}}
\newcommand{\bi}{\begin{itemize}}    \newcommand{\ei}{\end{itemize}}
\newcommand{\bt}{\begin{tabular}}    \newcommand{\et}{\end{tabular}}
\newcommand{\bdm}{\begin{displaymath}} \newcommand{\edm}{\end{displaymath}}

\newcommand{\D}{D\!\!\!\!/\,}

\newcommand{\ETSB}{\mathcal{S}\!\!\!\!/\,}
\newcommand{\nb}{\nabla\!\!\!\!/\,}

\def\qed{\hfill{q.e.d.}\smallskip\smallskip}

\begin{document}
	
	\title[]{Positive Energy Theorems for Spin Initial Data with Charge}       
	
	\author{Simon Raulot}
	\address[Simon Raulot]{Univ Rouen Normandie, CNRS, Normandie Univ, LMRS UMR 6085, F-76000 Rouen, France}
	\email{simon.raulot@univ-rouen.fr}
	
	\begin{abstract}
		We establish positive energy theorems for complete spin initial data sets with charge in dimensions $n \geq 4$, under a dominant energy condition and assuming the existence of at least one asymptotically flat end. Our results, formulated in the purely electric case, extend the classical theorems of Gibbons--Hull~\cite{GibbonsHull}, Gibbons--Hawking--Horowitz--Perry~\cite{GibbonsHawkingHorowitzPerry}, and Bartnik--Chru\'sciel~\cite{BartnikChrusciel}.
	\end{abstract}

	
	
	\thanks{}
	
	\date{\today}   
	
	\maketitle 
	\pagenumbering{arabic}

	
	\section{Introduction}
	
	
	The positive energy theorem is a foundational result in mathematical relativity. It deals with initial data sets for the Einstein equations and asserts that, under suitable asymptotic and energy conditions, the Arnowitt--Deser--Misner energy--momentum vector $(\E,\P) \in \mathbb{R}^{n,1}$, where $\E$ denotes the total energy and $\P \in \mathbb{R}^n$ the total linear momentum of the initial data set, is future--directed and causal in Minkowski space, that is $\E \geq |\P|$. Moreover, equality occurs if and only if the data arise from Minkowski spacetime. This theorem was first proved in dimension $3$ by Schoen and Yau \cite{SchoenYau1,SchoenYau6} using minimal surface techniques, and later by Witten \cite{Witten1} via a spinorial method based on the existence of solutions to a Dirac-type equation on a complete spin manifold.
	
	Witten’s approach extends naturally to higher dimensions, provided the manifold admits a spin structure, and has inspired a wide range of developments in geometric analysis and general relativity. In particular, it offers a flexible framework for incorporating additional matter fields, such as those arising in the Einstein--Maxwell theory.
	
	In the Einstein--Maxwell setting, the physical origins of the inequality trace back to Gibbons and Hull \cite{GibbonsHull} in the context of $N=2$ supergravity, and to Gibbons, Hawking, Horowitz, and Perry \cite{GibbonsHawkingHorowitzPerry} for spacetimes containing black holes. The presence of an inner boundary modeling an apparent horizon requires additional care by imposing natural boundary conditions, and Herzlich \cite{Herzlich5} and Bartnik and Chrusciel \cite{BartnikChrusciel} subsequently provided a rigorous analytic treatment of this result.
	
    In the time-symmetric setting, the charged positive energy theorem reduces to the estimate $m \geq |Q|$ for the ADM mass $m$ and the total charge $Q$, which follows from a variety of charged Penrose-type inequalities. In dimension three, Jaracz \cite{Jaracz} proved such an inequality for initial data with either an outermost minimal boundary or asymptotically cylindrical ends, via an adaptation of the inverse mean curvature flow method, while Bray, Hirsch, Kazaras and Khuri \cite{BrayHirschKazarasKhuri} obtained a related result using spacetime harmonic functions. In the non-extremal case with an outermost minimal surface boundary, Jang \cite{Jang} first proposed an argument for the charged Penrose inequality, which was later made rigorous by Huisken and Ilmanen \cite{HuiskenIlmanen} via the inverse mean curvature flow. The general case allowing disconnected horizons was subsequently established by Khuri, Weinstein, and Yamada \cite{KhuriWeinsteinYamada}, whose proof relies on the conformal flow introduced by Bray \cite{Bray}. In higher dimensions, de Lima, Gir\~ao, Loz\'orio, and Silva \cite{LimaGiraoLozorioSilva} derived Penrose-like inequalities for hypersurfaces isometrically embedded in Euclidean space.

	The equality case in the charged positive energy theorem is closely related to the notion of BPS states and is characterized by the existence of a super-covariantly constant spinor. Locally, Tod \cite{Tod83} classified all initial data admitting a super--covariantly constant spinor, showing that they arise from the Israel--Wilson--Perjès  class. Globally, Chruściel, Reall, and Tod \cite{ChruscielRealTod} proved that, under additional assumptions, these maximal solutions are precisely the Majumdar--Papapetrou spacetimes, including the extremal Reissner--Nordström solution. In a more general setting, Khuri and Weinstein \cite{KhuriWeinstein} investigated rigidity aspects by coupling a solution of the Dirac equation with the generalized Jang equation. 
	
	In this article, we establish a positive energy theorem for initial data sets with charge in arbitrary dimension $n \geq 4$, assuming the manifold is spin and contains an asymptotically flat end (see section \ref{MainDefinitions} for the precise definitions). We work in the purely electric case, where the magnetic field vanishes identically, and the data consist of a Riemannian manifold \((M^n, g)\), a symmetric tensor $K$, and a vector field $E$. Our main result is the following:
	\begin{theorem}\label{PositiveET}
		Let $n \geq 4$ and let $(M^n,g, K, E)$ be a complete spin initial data set with charge, containing at least one asymptotically flat end and satisfying the dominant energy condition (\ref{DEC}). Then $\E \geq \sqrt{|\P|^2 + Q^2}$. 
	\end{theorem}
	In particular, the ADM energy-momentum $(\E,\P)\in\mathbb{R}^{n,1}$ is causal and future-directed. In this situation, the {\it ADM mass} of the distinguished end is $m=\sqrt{\E^2-|\P|^2}$ and the theorem ensures $m\geq |Q|$. In the time-symmetric case $K=0$, mass and energy coincide, and the dominant energy condition reads 
	\begin{eqnarray}\label{DEC-TimeSymmetric}
		R\geq(n-1)(n-2)|E|^2+2(n-1)|{\rm div}(E)|,
	\end{eqnarray}
	where $R$ is the scalar curvature of the $(M^n,g)$, yielding the following positive mass theorem with charge:
	\begin{theorem}
		Let $n\geq 4$ and let $(M^n,g,E)$ be a complete time-symmetric initial data set (without boundary) with charge, containing at least one asymptotically flat end and satisfying the dominant energy condition (\ref{DEC-TimeSymmetric}) then $m\geq|Q|$. 
	\end{theorem}
	
	We also consider asymptotically flat initial data sets with charge carrying a compact inner boundary in which case we get a generalization, at least when the magnetic field is zero, of the well-known result of Gibbons, Hawking, Horowitz and Perry \cite{GibbonsHawkingHorowitzPerry} and  Bartnik and Chru\'sciel \cite{BartnikChrusciel}. In this situation, the strength of the gravitational field in the neighborhood of a $(n-1)$-hypersurface  $\Sigma$ may be measured by its null expansions 
	\begin{eqnarray}\label{FuturePastTrapped}
		\theta_\pm:=H\pm\Tr_\Sigma K 
	\end{eqnarray}
	where $H$ is the mean curvature with respect to the unit outward normal (pointing towards spatial infinity). The null expansions measure the rate of change of area for a shell of light emitted by the surface in the outward future direction (with $\theta_+$) and outward past direction (with $\theta_-$). Thus the gravitational field is interpreted as being strong near $\Sigma$ if $\theta_+\leq 0$ (resp. $\theta_-\leq 0$), in which
	case $\Sigma$ is referred to as a future (resp. past) trapped surface. Future or past 
	apparent horizons arise as boundaries of future or past trapped regions and satisfy the equation $\theta_+=0$ (resp. $\theta_-=0$). Then we prove:
	\begin{theorem}\label{PositiveET-Boundary}
		Under the assumptions of Theorem \ref{PositiveET}, suppose that $(M^n,g,K,E)$ has a compact inner boundary $\partial M$ whose connected components are either future or past trapped surfaces. Then the conclusion of Theorem \ref{PositiveET} still holds.
	\end{theorem}
	
	The key ingredient underlying all our arguments is a Schrödinger--Lichnerowicz formula adapted to initial data sets with charge in dimension $n \geq 4$, presented in Section \ref{SL-Section} (see Theorem \ref{SLforEM}). This approach is highly flexible and can be extended to various geometric and physical settings, such as Einstein--Maxwell data with negative cosmological constant or configurations with corners. These directions will be explored in future work. In Section \ref{Proofs}, we apply the analytical framework developed by Bartnik and Chruściel to establish our main results. This implies in particular that our positive energy theorems hold for manifolds that need only to have one asymptotically flat end and in particular there may be other complete ends but nothing about them is assumed other than the curvature condition.
	
	In this article we focus on establishing positivity. The equality case, however, requires a deeper analysis. In a forthcoming work~\cite{Raulot17}, we shall study in detail the geometry of solutions to $\overline{\nabla}\psi=0$ (see section \ref{ModifiedConnectionSection} for the precise definition), which we call \emph{charged generalized Killing spinors}. This approach will provide a general framework, independent of asymptotic conditions, and is expected to shed light on rigidity phenomena in the charged positive energy theorems.	
	
	
	\section{Notations and definitions}\label{MainDefinitions}
	
	
	An {\it initial data set with charge} $(M^n,g,K,E)$ is a Riemannian $n$-dimensional manifold $(M^n,g)$ e\-quip\-ped with a symmetric $(0,2)$-tensor $K$ and a vector field $E$ on $M$. We define the energy density $\mu$, current density $J$, and charge density $\varpi$ as
	\begin{equation}\label{EnergyDensity}
		\left\lbrace
		\begin{array}{cll}
			2\mu  & := & R-|K|^2+\Tr(K)^2-(n-1)(n-2)|E|^2 \\ 
			J & :=  & \div(K)-\nabla\Tr(K)\\
			\varpi & := & (n-1)\div(E)
		\end{array}
		\right.
	\end{equation}
	where ${\rm div}$ and $\nabla$ denote respectively the divergence operator on tensors and the Levi-Civita connection with respect to $g$. Then we say that $(M^n,g,K,E)$ satisfies the {\it dominant energy condition} if 
	\begin{eqnarray}\label{DEC}
		\mu\geq\sqrt{|J|^2+|\varpi|^2}.
	\end{eqnarray}
	We shall say that $M_{{\rm ext}}\subset M$ is an {\it asymptotically flat end} in the initial data set with charge $(M^n,g,K,E)$ if there exists a diffeomorphism $x:M_{{\rm ext}}\longrightarrow\mathbb{R}^n\setminus\overline{B}_1(0)$, where $\overline{B}_1(0)$ is the standard closed unit ball and such that if we denote $x=(x_1,...,x_n)$, the asymptotically coordinate chart, we have
	\begin{eqnarray}\label{AsymptoticDecay}
		g_{ij}(x)=\delta_{ij}+O_2(|x|^{-\tau}),\quad K_{ij}(x)=O_1(|x|^{-\tau-1})\quad\text{and}\quad E^i(x)=O_1(|x|^{-\tau-1}) 
	\end{eqnarray}
	for some $\tau>\frac{n-2}{2}$ and $1\leq i,j\leq n$. Moreover, we also require that $\mu$, $J$ and ${\rm div}(E)$ are integrable on $(M^n,g)$. Following Arnowitt, Deser and Misner \cite{ArnowittDeserMisner,ArnowittDeserMisner1}, we define the {\it ADM energy-momentum} $(\E,\P)\in\mathbb{R}^{n,1}$ of an end of $M$ by
	\begin{eqnarray*}
		\E & = & \frac{1}{2(n-1)\omega_{n-1}}\lim_{r\rightarrow\infty}\int_{S_r}\sum_{i,j=1}^n(g_{ij,i}-g_{ii,j})\overline{\nu}^jd\overline{\sigma}_r\\
		\P_i & = & \frac{1}{(n-1)\omega_{n-1}}\lim_{r\rightarrow\infty}\int_{S_r}\sum_{j=1}^n\big(k_{ij}-\Tr(K)g_{ij}\big)\overline{\nu}^jd\overline{\sigma}_r
	\end{eqnarray*}
	for $i=1,...,n$ where the right-hand sides are calculated in the asymptotically flat coordinate of the distinguished end and barred quantities are calculated using the Euclidean metric in the end. Moreover, $S_r$ denotes the Euclidean round sphere of radius $r>0$ and $\overline{\nu}$ its unit normal pointing toward infinity. Although this definition seems to depend on a specific choice of a coordinate chart, they are, as independently proven by Bartnik \cite{Bartnik1} and Chruściel \cite{Chrusciel1}, well-defined geometric invariants. In the presence of an electric field $E\in\Gamma(TM)$, the {\it total charge} is defined by
	\begin{eqnarray*}
		Q=\frac{1}{\omega_{n-1}}\lim_{r\rightarrow\infty}\int_{S_r}\sum_{j=1}^nE^j\overline{\nu}^jd\overline{\sigma}_r.
	\end{eqnarray*}
	
	A natural class of examples of such data arises from electrovacuum spacetimes that is $(n+1)$-dimensional time-oriented Lorentzian manifolds $(\mathcal{M}^{n+1},\mathfrak{g})$ satisfying the Einstein--Maxwell equations
	\begin{eqnarray*}\label{eq:EMfactor2}
		\mathfrak{Ric} - \frac{1}{2} \mathfrak{R}\,\mathfrak g =\mathcal{T}_F.
	\end{eqnarray*}
	Here $\mathfrak{Ric}$ and $\mathfrak{R}$ are respectively the Ricci tensor and the scalar curvatures of the spacetime and $\mathcal{T}_F$ is the energy-momentum tensor defined by
	\begin{eqnarray*}
		\mathcal{T}_F:= 2\big( F \circ F - \frac{1}{4} |F|^2 \mathfrak g\big),
	\end{eqnarray*} 
	where $F$ is a $2$-form, called the Faraday tensor, satisfying $dF = 0$ and $d\left(\ast F\right) = 0$. Moreover $(F\circ F)_{\alpha\beta}=\mathfrak g^{\sigma\gamma}F_{\alpha\sigma}F_{\beta\gamma}$ where Greek indices range from $1$ to $n+1$, and $\ast$ denotes the spacetime Hodge operator. Restricting to the purely electric case and assuming there exists an oriented spacelike hypersurface $M^n \subset \mathcal M$ with unit timelike normal $T$ and second fundamental form $K$, the Faraday tensor is expressed as $F = T^\flat \wedge E^\flat$, where $E$ is a vector field tangent to $M$. The induced initial data $(M^n,g,K,E)$ then satisfies the Einstein--Maxwell constraint equations:
	\begin{eqnarray*}
		R - |K|^2 + (\Tr K)^2 = (n-1)(n-2) |E|^2, \quad 
		\div(K)-\nabla\Tr(K)=0, \quad
		\div (E) = 0.
	\end{eqnarray*}
	
	Explicit solutions are provided by the Majumdar--Papapetrou family, describing multi-black-hole configurations in exact electrostatic equilibrium. Given $(x_j,m_j)\in\mathbb{R}^n\times\mathbb{R}^\ast_+$ for $j=1,...,l$, this spacetime has the manifold structure $\mathbb{R}\times\left(\mathbb{R}^n\setminus\{x_1,...,x_l\}\right)$ and is equipped with the metric
	\begin{eqnarray*}
		\mathfrak{g}_{MP} = -U^{-2} dt^2 + U^{2/(n-2)} \delta, \qquad U(x) = 1 + \sum_{j=1}^l \frac{m_j}{|x-x_j|^{n-2}},
	\end{eqnarray*}
	where $\delta$ is the $n$-dimensional Euclidean metric.  Here the Faraday tensor is given by $F=-dU^{-1}\wedge dt$ so that it is easily seen that mass and charge are equal. For a single center, that is $l=1$, this reduces to a spherically symmetric geometry corresponding to the extreme Reissner--Nordström solution. These spacetimes are extreme in the sense that they satisfy the equality case of our positive energy theorems.

	
	\section{The Schr\"odinger-Lichnerowicz formula for initial data sets with charge}\label{SL-Section}
	
	
	\subsection{Preliminaries on spinors}
	
	In this section, we consider $(M^n,g,K,E)$ an initial data set with charge endowed with a spin structure. Then there exists a smooth Hermitian vector bundle over $M$, the spinor bundle, denoted by $\SB$, whose sections are called spinor fields. The Hermitian scalar product is denoted by $\<\,,\,\>$. Moreover, the tangent bundle $TM$ (in fact the Clifford bundle) acts on $\SB$ by Clifford multiplication $X\otimes \psi\mapsto X\cdot\psi$ satisfying
	\begin{eqnarray*}\label{CliffordRule}
		X\cdot Y+Y\cdot X=-2g(X,Y)\,{\rm Id}
	\end{eqnarray*}
	and which is skew-Hermitian with respect to the Hermitian scalar product on $\SB$ that is 
	\begin{eqnarray*}\label{SkewHermitian}
		\<X\cdot\varphi,\psi\>  =  -\<\varphi,X\cdot\psi\>
	\end{eqnarray*}
	for any tangent vector fields $X$ and any spinor fields $\varphi$, $\psi\in\Gamma(\SB)$. On the other hand, the Riemannian Levi-Civita connection $\nabla$ lifts to the so-called spin Levi-Civita connection, also denoted by $\nabla$, and defines a metric covariant derivative on $\SB$ that preserves the Clifford multiplication. This means that 
	\begin{eqnarray*}\label{MetricCompatibility}
		X\<\varphi,\psi\>  =  \<\nabla_X\varphi,\psi\>+\<\varphi,\nabla_X\psi\>
	\end{eqnarray*}
	and 
	\begin{eqnarray*}\label{CliffordCompatibility}
		\nabla_X(Y\cdot\varphi)  =  \nabla_XY\cdot\varphi+Y\cdot\nabla_X\varphi
	\end{eqnarray*}
	hold for all $X$, $Y\in\Gamma(TM)$ and $\varphi$, $\psi\in\Gamma(\SB)$. A quadruplet $(\SB,\<\,,\,\>,\nabla,\,\cdot\,)$ satisfying these properties is called a Dirac bundle. The Dirac operator is then the first order elliptic differential operator acting on $\SB$ locally defined by
	\begin{eqnarray*}
		D\varphi=\sum_{j=1}^ne_j\cdot\nabla_{e_j}\varphi
	\end{eqnarray*}
	for $\varphi\in\Gamma(\SB)$. Here, and in all this work, $\{e_1,\cdots,e_n\}$ is a local orthonormal frame on $(M^n,g)$. 
	
	When $n$ is even, there exists a chirality operator, namely an endomorphism $\gamma$ of $\SB$ such that
	\begin{eqnarray}\label{ChiralityProperties}
		\<\gamma\varphi,\psi\>=\<\varphi,\gamma\psi\>,\quad\gamma^2={\rm Id},\quad \big\{ X,\gamma\big\}=0,\quad \nabla\gamma=0
	\end{eqnarray}
	where $\big\{ X,\gamma\big\}\varphi:=X\cdot\gamma\varphi+\gamma(X\cdot\varphi)$ for all $\varphi$, $\psi\in\Gamma(\SB)$ and $X\in\Gamma(TM)$. Such an operator is given by the Clifford action of the volume element of the spinor bundle $\SB$. When n is odd, a chirality operator may not exist.To overcome this, we consider $\SB\oplus\SB$, the direct sum of two copies of the spinor bundle over $M$, equipped with the Hermitian metric
	\begin{eqnarray*}
		\<(\varphi_1,\varphi_2),(\psi_1,\psi_2)\>:=\<\varphi_1,\psi_1\>+\<\varphi_2,\psi_2\>
	\end{eqnarray*}
	on which we define the Clifford action by
	\begin{eqnarray*}\label{TCM}
		X\bullet(\varphi_1,\varphi_2):=\left(X\cdot\varphi_1,-X\cdot\varphi_2\right)
	\end{eqnarray*}
	and the associated connection by $\nabla\oplus\nabla$. The endomorphism 
	\begin{eqnarray*}\label{ChiralityOdd}
		\gamma(\varphi_1,\varphi_2):=(\varphi_2,\varphi_1)
	\end{eqnarray*}
	then satisfies the same chirality properties (\ref{ChiralityProperties}) as in the even case. Thus by letting
	$$
	(\TSB,\<\,,\,\>,\nabla,\,\cdot\,):=
	\left\lbrace
	\begin{array}{ll}
		(\SB,\<\,,\,\>,\nabla,\,\cdot\,) & \text{ if } n \text{ is even}\\
		(\SB\oplus\SB,\<\,,\,\>\oplus\<\,,\,\>,\nabla\oplus\nabla,\,\bullet\,) & \text{ if } n \text{ is odd}
	\end{array}
	\right.
	$$
	we obtain a Dirac bundle on which there always exists $\gamma\in\Gamma({\rm End\,}\TSB)$ satisfying (\ref{ChiralityProperties}). In the following, $D$ denotes the associated Dirac operator which satisfies the well-known Schr\"odinger-Lichnerowicz formula 
	\begin{eqnarray}\label{SL-Formula}
		D^2\varphi=\nabla^*\nabla\varphi+\frac{R}{4}\varphi
	\end{eqnarray}
	for all $\varphi\in\Gamma(\TSB)$. Here $\nabla^*$ is the $L^2$-formal adjoint of the connection $\nabla$ and the rough Laplacian is locally given by
	\begin{eqnarray}\label{RL-Local}
		\nabla^*\nabla\varphi=-\sum_{j=1}^n \nabla_{e_j}\nabla_{e_j}\varphi.
	\end{eqnarray}
	For a bounded domain $\Omega\subset M$ with smooth boundary, we have 
	\begin{eqnarray}\label{Dirac-IPP}
		\int_\Omega\<D\varphi,\psi\>d\mu=\int_\Omega\<\varphi,D\psi\>d\mu+\int_{\partial\Omega}\<\nu\cdot\varphi,\psi\>d\sigma
	\end{eqnarray}
	where $\nu$ is the outward unit normal to $\partial\Omega$ and $d\mu$ (resp. $d\sigma$) is the Riemannian volume form of $\Omega$ (resp. $\partial\Omega$) with respect to $g$. 
	
	\subsection{The modified connection}\label{ModifiedConnectionSection}
	
	We define the following modified connection:
	\begin{eqnarray}\label{ModifiedConnection}
		\overline{\nabla}_X\varphi:=\nabla_X\varphi-\frac{1}{2}E\cdot X\cdot\gamma\varphi+\frac{n-3}{2}g(E,X)\gamma\varphi+\frac{1}{2}K(X)\cdot\gamma\varphi
	\end{eqnarray}
	for all $X\in \Gamma(TM)	$ and $\varphi\in\Gamma(\TSB)$ and the associated Dirac operator locally given by
	\begin{eqnarray}\label{ModifiedDirac}
		\overline{D}\varphi:=\sum_{j=1}^ne_j\cdot\overline{\nabla}_{e_j}\varphi.
	\end{eqnarray}
	\begin{remark}
		When $n=3$, the modified connection coincides with the one used by Gibbons, Hawking, Horowitz, and Perry \cite{GibbonsHawkingHorowitzPerry}, and later by Bartnik--Chru\'sciel \cite{BartnikChrusciel}, in the case where the magnetic field vanishes.
	\end{remark}
	
	We first begin by noticing the following properties of this operator.
	\begin{lemma}\label{MD-Properties}
		The operator $\overline{D}$ is an elliptic differential operator of order one, symmetric with respect to the $L^2$-scalar product on $\TSB$ which satisfies
		\begin{eqnarray}\label{RelationsDOD}
			\overline{D}\varphi=D\varphi-\frac{1}{2}\left(E+\Tr(K)\right)\cdot\gamma\varphi
		\end{eqnarray}
		for all $\varphi\in\Gamma(\TSB)$. 
	\end{lemma}
	
	{\it Proof.} The fact that $\overline{D}$ is an elliptic differential operator of order one follows directly from (\ref{RelationsDOD}) since then we observe that $\overline{D}$ is a zero order modification of the Dirac operator $D$. Now we compute that 
	\begin{eqnarray*}
		\overline{D} \varphi  = D\varphi-\frac{1}{2}\sum_{j=1}^n e_j\cdot E\cdot e_j\cdot\gamma\varphi
		+\frac{n-3}{2}\sum_{j=1}^n E^je_j\cdot\gamma\varphi+\frac{1}{2}\sum_{j=1}^ne_j\cdot K(e_j)\cdot\gamma\varphi
	\end{eqnarray*}
	where $E=\sum_{j=1}^nE^je_j$. It follows from the Clifford rule that
	\begin{eqnarray*}
		-\frac{1}{2}\sum_{j=1}^n e_j\cdot E\cdot e_j\cdot\gamma\varphi
		+\frac{n-3}{2}\sum_{j=1}^n E^je_j\cdot\gamma\varphi=-\frac{1}{2}E\cdot\gamma\varphi.
	\end{eqnarray*}
	On the other hand,  since $K$ is symmetric, we have
	\begin{eqnarray*}
		\frac{1}{2}\sum_{j=1}^ne_j\cdot K(e_j)\cdot\gamma\varphi=-\frac{1}{2}\Tr(K)\gamma\varphi.
	\end{eqnarray*}
	This leads to (\ref{RelationsDOD}). Finally, it is not difficult to check that the map 
	\begin{eqnarray*}
		\varphi\mapsto \big(E+\Tr(K))\cdot\gamma\varphi
	\end{eqnarray*}
	is pointwise symmetric so that the symmetry with respect to the $L^2$-scalar product follows from the one of the Dirac operator $D$. 
	\qed
	
	\begin{remark}\label{IPP-DM}
		When $\Omega$ is a bounded domain with smooth boundary $\partial\Omega$ in a spin initial data set with charge $(M^n,g,K,E)$, we deduce from (\ref{Dirac-IPP}) that 
		\begin{eqnarray*}
			\int_\Omega\<\overline{D}\varphi,\psi\>d\mu=\int_\Omega\<\varphi,\overline{D}\psi\>d\mu+\int_{\partial\Omega}\<\nu\cdot\varphi,\psi\>d\sigma
		\end{eqnarray*}
		for all $\varphi$, $\psi\in\Gamma(\TSB)$. 
	\end{remark}
	
	Now, in order to derive the Weitzenb\"ock formula associated to the modified Dirac operator $\overline{D}$, we need to express the $L^2$-formal adjoint of $\overline{\nabla}$. This is done in the following lemma. 
	\begin{lemma}\label{MC-Properties}
		The $L^2$-formal adjoint $\overline{\nabla}^\ast$ of the modified connection $\overline{\nabla}$ is given by
		\begin{eqnarray*}
			\overline{\nabla}^\ast\overline{\nabla}\varphi & = & \nabla^\ast\nabla\varphi-\frac{1}{2}\left\{D,E\cdot\gamma\right\}\varphi-\frac{n-1}{2}\div(E)\gamma\varphi-\frac{1}{2}\div(K)\cdot\gamma\varphi\\
			& & +\frac{1}{4}\big(|K|^2+\big((n-1)(n-2)+1\big)|E|^2\big))\varphi
		\end{eqnarray*}
		for all $\varphi\in\Gamma(\TSB)$.  
	\end{lemma}
	
	{\it Proof.} For $\varphi$, $\psi\in\Gamma(\TSB)$  we compute
	\begin{eqnarray*}
		\<\overline{\nabla}\varphi,\overline{\nabla}\psi\> & = & \div(\xi)+\underbrace{\sum_{j=1}^n\<-\nabla_{e_j}\overline{\nabla}_{e_j}\varphi,\psi\>}_{(1)} +\underbrace{\frac{1}{2}\sum_{j=1}^n\<-\overline{\nabla}_{e_j}\varphi,E\cdot e_j\cdot\gamma\psi\>}_{(2)}\\
		& & +\underbrace{\frac{n-3}{2}\sum_{j=1}^nE^j\<\overline{\nabla}_{e_j}\varphi,\gamma\psi\>}_{(3)}
		+\underbrace{\frac{1}{2}\sum_{j=1}^n\<\overline{\nabla}_{e_j}\varphi,K(e_j)\cdot\gamma\psi\>}_{(4)}
	\end{eqnarray*}
	where $\xi\in\Gamma(TM)$ is the vector field on $M$ defined by $g(\xi,X)=\<\overline{\nabla}_X\varphi,\psi\>$ for all $X\in\Gamma(TM)$. The first term in the previous identity can be written as
	\begin{eqnarray*}
		(1) & = & \underbrace{\sum_{j=1}^n\<-\nabla_{e_j}\nabla_{e_j}\varphi,\psi\>}_{(1a)}+\underbrace{\frac{1}{2}\sum_{j=1}^n\<\nabla_{e_j}\big(E\cdot e_j\cdot\gamma\varphi\big),\psi\>}_{(1b)}+\underbrace{\frac{n-3}{2}\sum_{j=1}^n\<-\nabla_{e_j}(E^j\gamma\varphi),\psi\>}_{(1c)}\\
		& &+\underbrace{\frac{1}{2}\sum_{j=1}^n\<-\nabla_{e_j}\big(K(e_j)\cdot\gamma\varphi\big),\psi\>}_{(1d)}.
	\end{eqnarray*}
	It follows from (\ref{RL-Local}) that $(1a)=\<\nabla^*\nabla\varphi,\psi\>$. On the other hand, we also have:
	\begin{eqnarray*}
		(1b) & = & -\frac{1}{2}\sum_{j=1}^n\<\nabla_{e_j}\big(e_j\cdot E \cdot\gamma\varphi+2E^j\gamma\varphi\big),\psi\>\\
		&  = & -\frac{1}{2}\sum_{j=1}^n\big(\<e_j\cdot\nabla_{e_j}(E\cdot\gamma\varphi)+2e_j(E^j)\gamma\varphi+2E^j\nabla_{e_j}(\gamma\varphi),\psi\>\big)\\
		& = & \<-\frac{1}{2}D(E\cdot\gamma\varphi)-\div(E)\gamma\varphi-\nabla_{E}(\gamma\varphi),\psi\>.
	\end{eqnarray*}
	Then using (\ref{CliffordCompatibility}) we compute
	\begin{eqnarray*}
		(1c) & = & -\frac{n-3}{2}\sum_{j=1}^n\<e_j(E^j)\gamma\varphi+E^j\nabla_{e_j}(\gamma\varphi),\psi\>\\
		& = & \<-\frac{n-3}{2}\div(E)\gamma\varphi-\frac{n-3}{2}\nabla_E(\gamma\varphi),\psi\>.
	\end{eqnarray*}
	Finally it holds that 
	\begin{eqnarray*}
		(1d) & = & -\frac{1}{2}\sum_{j=1}^n\<\nabla_{e_j}\big(K(e_j)\big)\cdot\gamma\varphi-\gamma\big(K(e_j)\cdot\nabla_{e_j}\varphi\big),\psi\>\\
		& = & \<-\frac{1}{2}\div(K)\cdot\gamma\varphi+\frac{1}{2}\gamma(D_K\varphi),\psi\>
	\end{eqnarray*}
	where we let $D_K:=\sum_{j=1}^nK(e_j)\cdot\nabla_{e_j}$. Putting all these identities together leads to 
	\begin{eqnarray*}
		(1)& = & \<\nabla^*\nabla\varphi-\frac{1}{2}D(E\cdot\gamma\varphi)-\frac{n-1}{2}\div(E)\gamma\varphi-\frac{1}{2}\div(K)\cdot\gamma\varphi,\psi\>\\
		&&+\<\frac{1}{2}\gamma(D_K\varphi)-\frac{n-1}{2}\nabla_E(\gamma\varphi),\psi\>
	\end{eqnarray*}
	Let us now examine the term $(2)$ which can be written as
	\begin{eqnarray*}
		(2) & = & \underbrace{\frac{1}{2}\sum_{j=1}^n\<-\nabla_{e_j}\varphi,E\cdot e_j\cdot\gamma\psi\>}_{(2a)}+\underbrace{\frac{1}{4}\sum_{j=1}^n\<E\cdot e_j\cdot\gamma\varphi,E\cdot e_j\cdot\gamma\psi\>}_{(2b)}\\
		& & +\underbrace{\frac{n-3}{4}\sum_{j=1}^nE^j\<-\gamma\varphi,E\cdot e_j\cdot\gamma\psi\>}_{(2c)}+\underbrace{\frac{1}{4}\sum_{j=1}^n\<-K(e_j)\cdot\gamma\varphi,E\cdot e_j\cdot\gamma\psi\>}_{(2d)}.
	\end{eqnarray*}
	Therefore we calculate
	\begin{eqnarray*}
		(2a) & = & -\frac{1}{2}\sum_{j=1}^n\<e_j\cdot E\cdot\nabla_{e_j}\varphi,\gamma\psi\>\\
		& = & \frac{1}{2}\sum_{j=1}^n\<E\cdot e_j\cdot\nabla_{e_j}\varphi+2E^j\nabla_{e_j}\varphi,\gamma\psi\>\\
		& = & \<-\frac{1}{2}E\cdot\gamma(D\varphi)+\nabla_E(\gamma\varphi),\psi\>. 
	\end{eqnarray*}
	On the other hand, it is simple to observe that $(2b)=(n/4)|E|^2\<\varphi,\psi\>$ and that
	\begin{eqnarray*}
		(2c)=-\frac{n-3}{4}\<\gamma\varphi,E\cdot E\cdot\gamma\psi\>=\frac{n-3}{4}|E|^2\<\varphi,\psi\>.
	\end{eqnarray*}
	Finally from the symmetry of $K$ and using the Clifford rule and (\ref{ChiralityProperties}), we easily check that
	\begin{eqnarray*}
		(2d) & = & -\frac{1}{4}\sum_{j=1}^n\<e_j\cdot E\cdot K(e_j)\cdot\gamma\varphi,\gamma\psi\>\\
		& = & \frac{1}{4}\sum_{j=1}^n\<E\cdot e_j\cdot K(e_j)\cdot\gamma\varphi+2E^jK(e_j)\cdot\gamma\varphi,\gamma\psi\>\\
		& = & \<\frac{1}{4}\Tr(K)E\cdot\varphi-\frac{1}{2}K(E)\cdot\varphi,\psi\>.
	\end{eqnarray*}
	Combining the fourth previous identities proves that
	\begin{eqnarray*}
		(2)=\<-\frac{1}{2}E\cdot\gamma(D\varphi)+\frac{2n-3}{4}|E|^2+\nabla_E(\gamma\varphi)+\frac{1}{4}\Tr(K)E\cdot\varphi-\frac{1}{2}K(E)\cdot\varphi,\psi\>.
	\end{eqnarray*}
	Let us now tackle the third term. For this, we write
	\begin{eqnarray*}
		(3) & = & \frac{n-3}{2}\<\overline{\nabla}_E\varphi,\gamma\psi\>\\
		& = &  \<\frac{n-3}{2}\nabla_E\varphi+\frac{n-3}{4}|E|^2\gamma\varphi+\frac{(n-3)^2}{4}|E|^2\gamma\varphi+\frac{n-3}{4}K(E)\cdot\gamma\varphi,\gamma\psi\>
	\end{eqnarray*}
	that is
	\begin{eqnarray*}
		(3)=\<\frac{n-3}{2}\nabla_E(\gamma\varphi)+\frac{(n-3)(n-2)}{4}|E|^2\varphi-\frac{n-3}{4}K(E)\cdot\varphi,\psi\>.
	\end{eqnarray*} 
	It remains to look at the fourth term which can be decomposed as follow:
	\begin{eqnarray*}
		(4) & = & \underbrace{\frac{1}{2}\sum_{j=1}^n\<\nabla_{e_j}\varphi,K(e_j)\cdot\gamma\psi\>}_{(4a)}+\underbrace{\frac{1}{4}\sum_{j=1}^n\<-E\cdot e_j\cdot\gamma\varphi,K(e_j)\cdot\gamma\psi\>}_{(4b)}\\
		& & +\underbrace{\frac{n-3}{4}\sum_{j=1}^nE^j\<\gamma\varphi,K(e_j)\cdot\gamma\psi\>}_{(4c)}+\underbrace{\frac{1}{4}\sum_{j=1}^n\<K(e_j)\cdot\gamma\varphi,K(e_j)\cdot\gamma\psi\>}_{(4d)}.
	\end{eqnarray*}
	It is direct to see that $(4a)=-(1/2)\<\gamma(D_K\varphi),\psi\>$. Moreover, from the computation of $(2d)$, we deduce that
	\begin{eqnarray*}
		(4b)=\<-\frac{1}{4}\Tr(K)E\cdot\varphi+\frac{1}{2}K(E)\cdot\varphi,\psi\>.
	\end{eqnarray*}
	Finally it is straightforward to observe that 
	\begin{eqnarray*}
		(4c) = \frac{n-3}{4}\<\gamma\varphi,K(E)\cdot\gamma\psi\>=\<\frac{n-3}{4}K(E)\cdot\varphi,\psi\>
	\end{eqnarray*}
	and $(4d)=(1/4)|K|^2\<\varphi,\psi\>$ so that
	\begin{eqnarray*}
		(4)=\<-\frac{1}{2}\gamma(D_K\varphi)-\frac{1}{4}\Tr(K)E\cdot\varphi+\frac{1}{2}K(E)\cdot\varphi+\frac{n-1}{4}K(E)\cdot\varphi+\frac{1}{4}|K|^2\varphi,\psi\>.
	\end{eqnarray*}
	The announced formula follows from the combination of the expressions of $(1)$, $(2)$, $(3)$ and $(4)$. 
	\qed
	
	\begin{remark}\label{IPP-Nabla}
		In the preceding proof, we proved the following pointwise equality:
		\begin{eqnarray*}
			\<\overline{\nabla}^*\overline{\nabla}\varphi,\psi\>= \<\overline{\nabla}\varphi,\overline{\nabla}\psi\>-\div(\xi)
		\end{eqnarray*}
		where $\xi\in\Gamma(TM)$ is the vector field defined by $g(\xi,X)=\<\overline{\nabla}_X\varphi,\psi\>$ for all $X\in\Gamma(TM)$ and where $\varphi$, $\psi\in\Gamma(\TSB)$. Then if $\Omega$ is a bounded domain with smooth boundary in a complete spin initial data set $(M^n,g,K,E)$, one can apply the divergence formula to obtain the following integration by parts formula:
		\begin{eqnarray*}
			\int_\Omega\<\overline{\nabla}^*\overline{\nabla}\varphi,\psi\>d\mu=\int_\Omega\<\overline{\nabla}\varphi,\overline{\nabla}\psi\>d\mu-\int_{\partial\Omega}\<\overline{\nabla}_\nu\varphi,\psi\>d\sigma.
		\end{eqnarray*}
	\end{remark}
	
	\subsection{The charged Schr\"odinger-Lichnerowicz formula}
	
	We can now state and prove the main formula of this section.
	\begin{theorem}\label{SLforEM}
		Let $(M^n,g,K,E)$ be a spin initial data set with charge, then
		\begin{eqnarray*}
			\overline{D}^2\varphi=\overline{\nabla}^*\overline{\nabla}\varphi+\mathcal{R}\varphi
		\end{eqnarray*}
		for any $\varphi\in\Gamma(\TSB)$ and where $\mathcal{R}\in\Gamma({\rm End\,}\TSB)$ is defined by
		\begin{eqnarray*}
			\mathcal{R}\varphi=\frac{1}{2}\big(\mu\varphi+\varpi\gamma\varphi+J\cdot\gamma\varphi\big).
		\end{eqnarray*}
	\end{theorem}
	
	{\it Proof.} To simplify the calculations, we define
	\begin{eqnarray*}
		\widetilde{D}\varphi:=D\varphi-\frac{1}{2}E\cdot\gamma\varphi
	\end{eqnarray*}
	for $\varphi\in\Gamma(\TSB)$. The modified Dirac operator can then be written as
	\begin{eqnarray*}
		\overline{D}\varphi= \widetilde{D}\varphi-\frac{1}{2}\Tr(K)\gamma\varphi.
	\end{eqnarray*}
	Then we compute 
	\begin{eqnarray}\label{FirstStep}
		\overline{D}^2\varphi=\widetilde{D}^2\varphi-\frac{1}{2}\Tr(K)\gamma(\widetilde{D}\varphi)-\frac{1}{2}\widetilde{D}\big(\Tr(K)\gamma\varphi\big)+\frac{1}{4}\Tr(K)^2\varphi.
	\end{eqnarray}
	However we have
	\begin{eqnarray*}
		-\frac{1}{2}\Tr(K)\gamma(\widetilde{D}\varphi) =-\frac{1}{2}\Tr(K)\gamma(D\varphi)-\frac{1}{4}\Tr(K)E\cdot\varphi
	\end{eqnarray*}
	and 
	\begin{eqnarray*}
		-\frac{1}{2}\widetilde{D}\big(\Tr(K)\gamma\varphi\big)= -\frac{1}{2}\nabla\Tr(K)\cdot\gamma\varphi+\frac{1}{2}\Tr(K)\gamma(D\varphi)+\frac{1}{4}\Tr(K)E\cdot\varphi.
	\end{eqnarray*}
	Putting these identities in (\ref{FirstStep}) leads to 
	\begin{eqnarray}\label{StepTwo}
		\overline{D}^2\varphi=\widetilde{D}^2\varphi -\frac{1}{2}\nabla\Tr(K)\cdot\gamma\varphi+\frac{1}{4}\Tr(K)^2\varphi.
	\end{eqnarray}
	A straightforward computation shows that
	\begin{eqnarray*}
		\widetilde{D}^2\varphi=D^2\varphi-\frac{1}{2}\big\{D,E\cdot\gamma\big\}\varphi+\frac{1}{4}|E|^2\varphi
	\end{eqnarray*}
	which, with the Schr\"odinger-Lichnerowicz formula (\ref{SL-Formula}), leads to
	\begin{eqnarray}\label{StepThree}
		\widetilde{D}^2\varphi=\nabla^*\nabla\varphi-\frac{1}{2}\big\{D,E\cdot\gamma\big\}\varphi+\frac{1}{4}(R+|E|^2)\varphi.
	\end{eqnarray}
	Combining (\ref{StepTwo}) and (\ref{StepThree}) and using Lemma \ref{MC-Properties} together with the definitions (\ref{EnergyDensity}) yields the claimed formula.
	\qed
	
	As a direct corollary of this formula, we obtain the following integral version for compact domains.
	\begin{corollary}\label{IntegralVersion}
		Let $\Omega$ be a bounded open set with smooth boundary in a complete spin initial data with charge $(M^n,g,K,E)$. Then for any $\varphi\in\Gamma(\TSB)$,
		\begin{eqnarray*}
			\int_\Omega\Big(|\overline{\nabla}\varphi|^2+\<\mathcal{R}\varphi,\varphi\>-|\overline{D}\varphi|^2\Big)d\mu = \int_{\partial\Omega}\<\overline{L}_{\nu}\varphi,\varphi\>d\sigma
		\end{eqnarray*}
		where $\nu$ is the outward unit normal to $\partial\Omega$ and $\overline{L}_\nu\varphi:=\overline{\nabla}_{\nu}\varphi+\nu\cdot\overline{D}\varphi$.
	\end{corollary}
	
	{\it Proof.} It is sufficient to integrate the formula from Theorem \ref{SLforEM} over the domain $\Omega$, and then use Remarks \ref{IPP-DM} and \ref{IPP-Nabla}.
	\qed
	
	
	\section{Proof of the positive energy theorems}\label{Proofs}
	
	
	We are now ready to prove the positive energy theorems using Witten's approach. The first step is to write the ADM energy-momentum $(\E,\P)$ and the total charge $Q$ of the distinguished end in terms of spinors. 
	\begin{proposition}\label{EMC-Spinors}
		Let $(M^n,g,K,E)$ be an initial data set with charge containing a distinguished asymptotically flat end, and let $e_1,...,e_n$ be an orthonormal frame near infinity. There exists $\psi_0\in\Gamma(\TSB)$ which is constant with respect to this frame such that 
		\begin{eqnarray*}
			\lim_{r\rightarrow\infty}\int_{S_r}\<\overline{L}_{\nu_r}\psi_0,\psi_0\>d\sigma_r=\frac{n-1}{2}\omega_{n-1}\big(\E-\sqrt{|\P|^2+Q^2}\big)
		\end{eqnarray*}
		where $\nu_r$ is the outward unit normal to $S_r$ and $\overline{L}_{\nu_r}$ is defined in Corollary \ref{IntegralVersion}. 
	\end{proposition} 
	
	{\it Proof.} We first choose $\psi_0\in\Gamma(\TSB)$ any spinor which is constant with respect to the chosen frame. Then we compute
	\begin{eqnarray*}
		\overline{L}_{\nu}\psi_0 & = & \sum_{i,j=1}^n\nu^i(\delta_{ij}+e_i\cdot e_j)\cdot\overline{\nabla}_{e_j}\psi_0\\
		& = & L_\nu\psi_0+\frac{1}{2}\sum_{i,j=1}^n\big(K_{ij}-\Tr(K)\delta_{ij}\big)e_j\cdot\gamma\psi_0+\frac{n-1}{2}g(E,\nu)\gamma\psi_0
	\end{eqnarray*}
	and $L_\nu\psi_0=\nabla_\nu\psi_0+\nu\cdot D\psi_0$. Then it is well-known (see \cite{Lee} for example) that under our decay assumptions it holds that
	\begin{eqnarray*}
		\lim_{r\rightarrow\infty}\int_{S_r}\<L_\nu\psi_0,\psi_0\>d\sigma_r  =  \frac{n-1}{2}\omega_{n-1}\E|\psi_0|^2
	\end{eqnarray*}
	and 
	\begin{eqnarray*}
		\lim_{r\rightarrow\infty}\int_{S_r}\<\frac{1}{2}\sum_{i,j=1}^n\big(K_{ij}-\Tr(K)\delta_{ij}\big)e_j\cdot\gamma\psi_0,\psi_0\>d\sigma_r  =   \frac{n-1}{2}\omega_{n-1}\<\sum_{j=1}^n \P_je_j\cdot\gamma\psi_0,\psi_0\>.
	\end{eqnarray*}
	On the other hand, it is straightforward to see that 
	\begin{eqnarray*}
		\lim_{r\rightarrow\infty}\int_{S_r}g(E,\nu)\<\gamma\psi_0,\psi_0\>d\sigma_r=\omega_{n-1}Q\<\gamma\psi_0,\psi_0\>.
	\end{eqnarray*}
	Putting all together yields the desired result. 
	\qed
	
	The second step of the proof consists to show that the modified Dirac operator is an isomorphism between a certain Hilbert space $\mathbb{H}$ and $L^2$, the space of square integrable sections of the bundle $\TSB$. For this, we adopt the framework developed by Bartnik and Chru\'sciel \cite{BartnikChrusciel}.  We first observe that under our assumptions, we have a {\it weighted Poincar\'e inequality} for the connection $\overline{\nabla}$. More precisely, this means that there exists $w\in L^1_{loc}$ with ${\rm ess\,inf}_\Omega w>0$ for all relatively compact $\Omega$ in $M$ such that for all $\varphi\in C^1_c$, the space of compactly supported $C^1$ spinor fields on $M$, we have 
	\begin{eqnarray}\label{wPi}
		\int_M|\varphi|^2w\,d\mu\leq\int_M |\overline{\nabla}\varphi|^2d\mu.
	\end{eqnarray}
	Indeed, since $\overline{\Gamma}_S$, the symmetric part of the connection $\overline{\nabla}$, is given by 
	\begin{eqnarray*}
		\overline{\Gamma}_S(X)=-\frac{1}{2}\big(K(X)+(n-2)g(E,X)\big)\gamma\in\Gamma\big({\rm End\,}\TSB\big)
	\end{eqnarray*}
	it is straightforward to check that the assumptions of \cite[Theorem 9.10]{BartnikChrusciel} are fulfilled because of the decay conditions (\ref{AsymptoticDecay}). On the other hand, if we assume that the dominant energy condition (\ref{DEC}) is fulfilled, the curvature endomorphism $\mathcal{R}$ appearing in Theorem \ref{SLforEM} is nonnegative since 
	\begin{eqnarray}\label{DEC1}
		\<\mathcal{R}\varphi,\varphi\>=\frac{1}{2}\<\big(\mu\varphi+\varpi\gamma\varphi+J\cdot\gamma\varphi\big),\varphi\>\geq \frac{1}{2}\big(\mu-\sqrt{|J|^2+|\varpi|^2}\big)|\varphi|^2\geq 0
	\end{eqnarray}
	for all $\varphi\in\Gamma(\TSB)$.
	
	\subsection{The boundaryless case} 
	
	We will first assume that $M$ has no compact inner boundary. In this situation, we have the following {\it Schr\"odinger-Lichnerowicz estimate} for the Dirac operator $\overline{D}$ (in the sense of \cite{BartnikChrusciel}) that is
	\begin{eqnarray}\label{SL-estimate}
		\int_M|\overline{\nabla}\varphi|^2d\mu\leq \int_M|\overline{D}\varphi|^2d\mu
	\end{eqnarray}
	for all $\varphi\in C^1_c$. This follows directly by integrating the formula of Theorem \ref{SLforEM} and by (\ref{DEC1}). This last property also implies that 
	\begin{eqnarray*}
		||\varphi||^2 :=\int_{M}\big(|\overline{\nabla}\varphi|^2+\<\mathcal{R}\varphi,\varphi\>\big)d\mu
	\end{eqnarray*}
	defines a norm on $C^1_c$. Therefore the space
	\begin{eqnarray*}
		\mathbb{H}:=||\,.\,||-{\rm completion\,\,of\,\,}C^1_c
	\end{eqnarray*}
	is a Hilbert space. The Poincar\'e inequality (\ref{wPi}) ensures that $\mathbb{H}$ embeds continuously in $H^1_{loc}$. In particular, it implies that any $\varphi\in\mathbb{H}$ can be represented by a spinor field in $H^1_{loc}$. Consider now the bilinear form defined by
	\begin{eqnarray*}
		\alpha(\varphi,\psi):=\int_M\<\overline{D}\varphi,\overline{D}\psi\>d\mu
	\end{eqnarray*}
	for $\varphi$, $\psi\in\mathbb{H}$. From Lemma $8.5$ in \cite{BartnikChrusciel}, we get that the map $\varphi\in\mathbb{H}\mapsto \overline{D}\varphi\in L^2$ is continuous and so $\alpha$ is also continuous on $\mathbb{H}\times\mathbb{H}$. Moreover, using the weighted Poincar\'e inequality (\ref{wPi}) and the Schr\"odinger-Lichnerowicz estimate (\ref{SL-estimate}), we immediately observe that $\alpha$ is coercive on $\mathbb{H}$. So if for $\chi\in L^2$, we define the continuous linear form
	\begin{eqnarray*}
		F_\chi(\varphi)=\int_M\<\chi,\overline{D}\varphi\>d\mu
	\end{eqnarray*}
	on $\mathbb{H}$, the Lax-Milgram theorem implies that there exists a unique $\xi_0\in\mathbb{H}$ such that $F_\chi(\varphi)=\alpha(\xi_0,\varphi)$ for all $\varphi\in\mathbb{H}$. In other words, we get that $\zeta:=\overline{D}\xi_0-\chi\in L^2$ is a weak solution of $\overline{D}\zeta=0$ since, from Remark \ref{IPP-DM}, the Dirac-type operator $\overline{D}$ is symmetric with respect to the $L^2$-scalar product. From the ellipticity of $\overline{D}$, we conclude that $\zeta\in\mathbb{H}\cap L^2$ is in fact a strong solution of this equation. Now it follows from standard arguments that, under the assumptions of Theorem \ref{PositiveET}, the operator $\overline{D}$ has a trivial $L^2$-kernel, so $\zeta\equiv 0$ and $\xi_0\in\mathbb{H}$ is the unique solution of the Dirac equation $\overline{D}\xi_0= \chi$. To summarize, we proved the following result:
	\begin{proposition}\label{DiracIsomorphismWB}
		Under the assumptions of Theorem \ref{PositiveET}, the operator $\overline{D}:\mathbb{H}\rightarrow L^2$ is an isomorphism.
	\end{proposition}
	Now we can apply the classical Witten's argument to conclude. Take a constant spinor $\psi_0$ as in Proposition \ref{EMC-Spinors} and extend it as a smooth spinor field on $M$ with support in the distinguished end. From our asymptotic assumptions (\ref{AsymptoticDecay}), we observe that $\overline{D}\psi_0$ is $L^2$ on $(M^n,g)$. Then it follows from Proposition \ref{DiracIsomorphismWB} that there exists an unique $\varphi\in\mathbb{H}$ such that $\overline{D}\varphi=-\overline{D}\psi_0$. In other words, the spinor field $\psi=\varphi+\psi_0$ is $\overline{D}$-harmonic.
	If $\varphi$ was an element of $C^1_c$ one could directly apply Corollary \ref{IntegralVersion} and Proposition \ref{EMC-Spinors} to conclude that
	\begin{eqnarray}\label{mass-formula}
		\frac{n-1}{2}\omega_{n-1}\big(\E-\sqrt{|\P|^2+Q^2}\big)=\int_M\Big(|\overline{\nabla}\psi|^2+\<\mathcal{R}\psi,\psi\>\Big)d\mu\geq 0
	\end{eqnarray}
	because of (\ref{DEC1}). Actually, one can then show that the previous equality holds for $\varphi\in\mathbb{H}$ since the right-hand side of (\ref{mass-formula}) is continuous on $\mathbb{H}$ (and $C^1_c$ is dense in $\mathbb{H}$). This concludes the proof of Theorem \ref{PositiveET}.
	
	\subsection{Positive energy theorem for charged black holes} 
	
	In this setting, we will apply the integral formula of Corollary \ref{IntegralVersion} to the bounded domain $\Omega_r$ of $M$ whose boundary is the union of the inner boundary $\partial M$ and large spheres $S_r$. To control the additional boundary term coming from the inner boundary, we need to introduce some notations. On the restricted spinor bundle $\ETSB:=\TSB_{|\partial M}$, we define for $X\in\Gamma(T\partial M)$ and $\varphi\in\Gamma(\ETSB)$ the linear connection 
	\begin{eqnarray*}
		\nb_X\varphi :=\nabla_X\varphi+\frac{1}{2} A(X)\cdot\nu\cdot\varphi
	\end{eqnarray*}
	where $A(X):=\nabla_X\nu$ is the Weingarten map of $\partial M$ in $M$ and the associated Dirac operator, denoted by $\D$, is then as usual locally given by
	\begin{eqnarray*}
		\D\varphi=\sum_{i=1}^{n-1}e_i\cdot\nu\cdot\nb_{e_i}\varphi
	\end{eqnarray*}
	for all $\varphi\in\Gamma(\ETSB)$. Here $\nu$ denotes the unit normal to $\partial M$ pointing toward infinity. We can therefore express the inner boundary term as follow:
	\begin{proposition}
		For all $\varphi\in\Gamma(\ETSB)$, the following identity holds
		\begin{eqnarray*}
			\overline{L}_\nu\varphi=-\D\varphi-\frac{1}{2}\mathcal{H}\varphi
		\end{eqnarray*}
		where $\mathcal{H}\in\Gamma({\rm End\,}\ETSB)$ is defined by
		\begin{eqnarray*}
			\mathcal{H}\varphi:=H\varphi+\Tr_{\partial M}(K)\nu\cdot\gamma\varphi-K(\nu)^\top\cdot\gamma\varphi-(n-1)E_\nu\gamma\varphi
		\end{eqnarray*}
		where $H$ is the mean curvature of $\partial M$, $\Tr_{\partial M}(K)$ is the trace along $\partial M$ of $K$, $K(\nu)^\top$ is the tangent part of the vector field $K(\nu)$ defined along $\partial M$ by $g\big(K(\nu),X\big)=K(X,\nu)$ for $X\in\Gamma(TM_{|\partial M})$ and $E_\nu=g(E,\nu)$.
	\end{proposition}
	{\it Proof.} Using formula (\ref{ModifiedConnection}) for the modified connection together with identity (\ref{RelationsDOD}), one obtains
	\begin{eqnarray*}
		\overline{L}_\nu\varphi   = 	\nabla_\nu\varphi +\nu\cdot D\varphi-\frac{1}{2}\{\nu,E\}-\frac{1}{2}\Tr(K)\nu\cdot\gamma\varphi+\frac{n-3}{2}E_\nu\gamma\varphi+\frac{1}{2}K(\nu)\cdot\gamma\varphi.	
	\end{eqnarray*}
	On the other hand, a standard computation shows that 
	\begin{eqnarray*}
		\nabla_\nu\varphi +\nu\cdot D\varphi=-\D\varphi-\frac{1}{2}H\varphi.
	\end{eqnarray*}
	Finally, using the Clifford rule together with the decompositions 
	\begin{eqnarray*}
		\Tr(K)=\Tr_{\partial M}(K)+K(\nu,\nu), 
		\qquad 
		K(\nu)=K(\nu)^\top+K(\nu,\nu)\nu,
	\end{eqnarray*}
	the claimed formula follows immediately.
	\qed 
	
	Using this expression in Corollary \ref{IntegralVersion}, we deduce the following integral version of Theorem \ref{SLforEM}:
	\begin{eqnarray}\label{IntegralVersionBoundary}
		\int_{S_r}\<\overline{L}_{\nu_r}\varphi,\varphi\>d\sigma_r = \int_{\Omega_r}\Big(|\overline{\nabla}\varphi|^2+\<\mathcal{R}\varphi,\varphi\>-|\overline{D}\varphi|^2\Big)d\mu - \int_{\partial M}\<\D\varphi+\frac{1}{2}\mathcal{H}\varphi,\varphi\>d\sigma
	\end{eqnarray}
	for all $\varphi\in\Gamma(\TSB)$. To proceed as in the boundaryless case, one needs to impose boundary conditions that make the modified Dirac operator $\overline{D}$ both elliptic and symmetric. For this purpose, we consider the pointwise projections
	\begin{eqnarray*}
		\pi_\pm \varphi := \frac{1}{2} (\varphi\pm\nu \cdot \gamma \varphi), \qquad \varphi \in \Gamma(\ETSB),
	\end{eqnarray*}
	which are well-known to define elliptic boundary conditions for the standard Dirac operator $D$ (see \cite{BartnikChrusciel} for example). Since ellipticity depends only on the principal symbol and $D$ and $\overline{D}$ share the same principal symbol (see Lemma \ref{MD-Properties}), these conditions are also elliptic for $\overline{D}$. Moreover, the symmetry of $\overline{D}$ under the boundary conditions $\pi_\pm$ follows directly from the Green’s formula established in Remark \ref{IPP-DM}. Now let 
	\begin{eqnarray*}
		\widetilde{C}^1_c:=\big\{\varphi\in C^1_c\,|\, \pi_-\varphi_{|\partial M_+}=0\text{ and }\pi_+\varphi_{|\partial M_-}=0\big\}
	\end{eqnarray*}
	where $\partial M =\partial M_+\cup\partial M_-$ with $\partial M_\pm$ is the portion of the boundary with $\theta_\pm\leq 0$. Here $\theta_\pm$ denote the null expansions of $\partial M$ defined by (\ref{FuturePastTrapped}). Then we prove:
	\begin{proposition}
		Under the assumptions of Theorem \ref{PositiveET-Boundary}, the Schr\"odinger-Lichnerowicz estimate (\ref{SL-estimate}) holds for all $\varphi\in\widetilde{C}_c^1$. 
	\end{proposition}
	{\it Proof.} A direct consequence of Formula (\ref{IntegralVersionBoundary}) and the dominant energy condition (\ref{DEC1}) is that 
	\begin{eqnarray*}
		\int_M|\overline{D}\varphi|^2d\mu\geq \int_{M}|\overline{\nabla}\varphi|^2d\mu - \int_{\partial M}\<\D\varphi+\frac{1}{2}\mathcal{H}\varphi,\varphi\>d\sigma
	\end{eqnarray*}
	for all $\varphi\in\widetilde{C}^1_c$. Now it is not difficult to check using the definition of $\D$ and (\ref{ChiralityProperties}) that
	\begin{eqnarray*}
		\<\D\varphi,\varphi\> = \<\D(\pi_+\varphi),\pi_-\varphi\>+\<\D(\pi_-\varphi),\pi_+\varphi\>
	\end{eqnarray*}
	for all $\varphi\in\Gamma(\TSB)$ so that $\<\D\varphi,\varphi\>=0$ as soon as $\varphi\in\widetilde{C}^1_c$. Similarly, it is straightforward to check that
	\begin{eqnarray*}
		\<\nu\cdot\gamma\varphi,\varphi\>  =  |\pi_+\varphi|^2-|\pi_-\varphi|^2,\quad\<\gamma\varphi,\varphi\> =  2 {\rm Re} \<\gamma(\pi_+\varphi),\pi_-\varphi\>
	\end{eqnarray*}
	and 
	\begin{eqnarray*} 
		\<X\cdot\gamma\varphi,\varphi\>  = 2{\rm Re} \<X\cdot\gamma(\pi_+\varphi),\pi_-\varphi\>
	\end{eqnarray*}
	for all $X\in\Gamma(T\partial M)$ and $\varphi\in\Gamma(\TSB)$. Here ${\rm Re} \<\,,\,\>$ denotes the real part of the Hermitian scalar product on $\TSB$. Thus it holds on $\partial M_\pm$ that $\<\mathcal{H}\varphi,\varphi\>=\theta_\pm|\pi_{\pm}\varphi|^2$ for all $\varphi\in\widetilde{C}^1_c$. Combining the previous estimates leads to 
	\begin{eqnarray*}
		\int_M|\overline{D}\varphi|^2d\mu\geq \int_{M}|\overline{\nabla}\varphi|^2d\mu -\frac{1}{2}\int_{\partial M_+}\theta_+|\pi_+\varphi|^2d\sigma-\frac{1}{2}\int_{\partial M_-}\theta_-|\pi_-\varphi|^2d\sigma
	\end{eqnarray*}
	and this concludes the proof since $\theta_\pm\leq 0$ on $\partial M_{\pm}$. 
	\qed
	
	The proof of Theorem \ref{PositiveET-Boundary} proceeds exactly as in the previous section by taking $\mathbb{H}$ to be the $||\,.\,||-$completion of $\widetilde{C}^1_c$.
	
	\section{Perspectives}
	
	The results established here provide a natural extension of Witten’s spinorial method 
	to charged initial data in arbitrary dimension. Beyond the positivity statement itself, 
	a central geometric feature is the existence of charged generalized Killing spinors, 
	arising from the modified connection. A systematic study of the geometry of manifolds 
	admitting such spinors, not necessarily with asymptotically flat ends, will be carried 
	out in future work, with particular emphasis on applications to rigidity questions for 
	the charged positive energy theorem. 
	
	We also note that analogous formulas can be derived in other settings of interest, 
	including initial data with negative cosmological constant and initial data with corners. 
	Since these involve additional analytic difficulties, we postpone their detailed analysis 
	to forthcoming work.


	\bibliographystyle{alpha}     
	\bibliography{BiblioHabilitation}


\end{document}